\documentclass[12pt, leqno]{article}
\usepackage{amssymb,amsmath}
\evensidemargin =\oddsidemargin

\font\teneufm=eufm10
\font\seveneufm=eufm7
\font\fiveeufm=eufm5
\newfam\eufmfam
\textfont\eufmfam=\teneufm
\scriptfont\eufmfam=\seveneufm
\scriptscriptfont\eufmfam=\fiveeufm

\let\goth\frak

\def\gg{\goth g}
\def\gs{\goth s}
\def\gh{\goth h}
\def\gk{\goth k}
\def\gl{\goth l}

\def\gc{\goth c}
\def\gb{\goth b}
\def\gp{\goth p}

\def\gn{\goth n}
\def\gu{\goth u}
\def\gq{\goth q}
\def\gz{\goth z}
\def\go{\goth o}
\def\ga{\goth a}
\def\gd{\goth d}

\def\pp{\mbox{\bf p}}

\def\beq{\begin{equation}}
\def\eeq{\end{equation}}
\def\bea{\begin{eqnarray}}
\def\eea{\end{eqnarray}}
\def\beas{\begin{eqnarray*}}
\def\eeas{\end{eqnarray*}}

\def\cplus{\hbox{$\supset${\raise1.05pt\hbox{\kern -0.55em
${\scriptscriptstyle +}$}}\ }}

\DeclareMathOperator{\supp}{supp}

\DeclareMathOperator{\tr}{tr}

\DeclareMathOperator{\ch}{ch}

\DeclareMathOperator{\Sup}{Sup}
\DeclareMathOperator{\Ind}{Ind}

\newtheorem{theorem}[equation]{Theorem}

\newtheorem{lemma}[equation]{Lemma}

\newtheorem{corollary}[equation]{Corollary}

\newtheorem{proposition}[equation]{Proposition}

\def\Z{\mathbb Z}

\begin{document}

\title{On the structure and characters of weight modules}
\author{Dimitar Grantcharov\thanks{Research partially supported by an NSF GIG Grant}}
\date{}

\maketitle
\begin{abstract}
Let $\mathfrak g$ be a classical Lie superalgebra of type I or a
Cartan-type Lie superalgebra {\bf W}$(n)$. We study weight
$\gg$-modules using  a  method inspired by  Mathieu's
 classification of the simple weight
modules with finite weight multiplicities  over reductive Lie
algebras, \cite{M}. Our approach is based on the fact that every
simple weight $\mathfrak g$-module with finite weight multiplicities is
obtained via a composition  of a twist and localization from a
highest weight module.  This allows us to transfer many results
for  category ${\cal O}$ modules to the  category of  weight modules
with finite weight multiplicities. As a  main application of the
method we reduce the problems of finding a ${\mathfrak g}_0$-composition
series and a character formula for all simple weight modules  to the same problems for simple  highest
weight modules. In this way, using results of Serganova  we obtain a character formula for all simple weight
{\bf W}$(n)$-modules and all simple atypical nonsingular ${\mathfrak s}{\mathfrak l} (m|1)$-modules. Some of our
results are new already in the case of a classical reductive Lie
algebra $\mathfrak g$.
\end{abstract}

Key words (2000 MSC): Primary 17B10, Secondary 17B66.

\section{Introduction}

Let $\gg$ be a classical Lie superalgebra of type I, or let $\gg =
${\bf W}$(n)$. The problem of finding the characters of all finite
dimensional simple $\gg$-modules $M$ has been known for more than
20 years to be very interesting and  challenging. Many
mathematicians, including J. Bernstein, V. Kac, D. Leites, I.
Penkov, and V. Serganova have contributed to its solution, until
it was completely solved by Serganova for $\gg = \gg \gl (m|n)$.
Recently the infinite dimensional highest weight modules over Lie
superalgebras have attracted attention, as new methods of
J. Brundan and V. Serganova (see \cite{B}, \cite{S1}, and
\cite{S2}) have opened new perspectives for the study of these
modules.

 All modules considered in the present  paper are weight modules, i.e. are isomorphic to the direct
  sum of their weight spaces. Our main goal is to study
  quotients of parabolically induced $\gg$-modules $M$, in
  particular all simple weight $\gg$-modules with finite weight multiplicities. Let $\ga$ be a reductive Lie subalgebra of $\gg_0$
 containing a fixed Cartan subsuperalgebra $\gh$ of $\gg$, such that all elements $ a \in \ga, a \notin \gh$, act injectively on $M$. We
  focus our attention on  the $\gd$-module structure of $M$, where $\gd = \gg_0$, $\gd =
  \ga$ or in the case when $M$ is not simple, $\gd = \gg$. We show
  that the $\gd$-module structure of $M$ is described up to a partially
  finite enlargement (see section 4 for the definition of a partially finite
  enlargement) by the $\gd$-module  structure of a quotient $M_h$ of a Verma
  module, i.e. a highest weight module. The correspondence $M \mapsto
  M_h$ is a basic tool in Mathieu's classification of all simple
  weight modules over reductive Lie algebras. This correspondence
  is a specific  version of M. Duflo's theorem stating that the annihilator of every
  simple module over a simple Lie algebra  is the annihilator of a certain highest
  weight module with the same central character. For Lie superalgebras $\gg$ of type
  I and modules with bounded weight multiplicities
  the correspondence $M \mapsto M_h$ is studied  in \cite{G} where the ``finite
  enlargement phenomenon'' occurs when $\gd = \ga = \gg_0$. The aim of the present paper is twofold - to study
  this phenomenon at a larger scale and to obtain an explicit character
  formula for $M$ provided that such a formula is known for
  $M_h$.

One of our main tools is the localization of weight modules
determined by a commuting set of roots, introduced in \cite{M}.
Section 2 is devoted to the basic properties of localization. The
results in this section hold for arbitrary weight modules, which
gives us a hope that eventually the method will provide new
results for weight modules with infinite weight multiplicities, in
particular for Harish-Chandra modules. In Section 3 we study
$(\gg,\ga)$-coherent families. Our results generalizes 
Mathieu's theory of coherent families for
reductive Lie algebras and its extension in \cite{G} (where $\ga =
\gg_0$) for Lie superalgebras of type I. We show that every simple
weight $\gg$-module with finite weight multiplicities $M $ has a
unique semisimple $(\gg,\ga)$-coherent extension ${\cal
E}_{\ga}(M)$. As in the Lie algebra case,  ${\cal E}_{\ga}(M)$
contains a $D$-highest weight submodule with respect to any basis
$D$ of the root system $\Delta$ of $\gg$. In the fourth section we
find a necessary and sufficient condition for $M$ to be
$\gd$-semisimple module and apply it to the particular case $\gg =
\go \gs \gp (2|n)$. Our main result is contained in the  last
section. We write down a character formula for every simple weight
$\gg$-module with finite weight multiplicities in terms of composition multiplicities of Verma
$\gg$-modules  in the category ${\cal O}$. For the case of a
reductive Lie algebra $\gg$ this result is a new more explicit
version of a  character formula of  Mathieu. In the general case, i.e. for Lie superalgebras $\gg$ of type I and $\gg =${\bf W}$(n)$, we combine our results with a set of results of Serganova in \cite{S1}, \cite{S2}, and \cite{S3}. As a corollary we establish a character formula for all simple weight $\gg$-modules
 except those which have singular atypical central character for the cases of $\gg = \gs \gl (m|1)$ and $\gg = \go \gs \gp(2|n)$. Finally, for $\gg = ${\bf W}$(n)$, a character formula for all simple weight $\gg$-modules is provided.

\medskip
\noindent
{\bf Acknowledgments:} I thank V. Serganova for the numerous stimulating discussions and helpful comments, and G. Zuckerman for a suggestion which led to the remark at the end of Section 3. Special thanks are due to I. Penkov for his  careful reading the manuscript and suggesting valuable improvements.

\section{Localization of Weight Modules}

\subsection{Notations and Conventions} The ground field $K$ is  algebraically
 closed and of characteristic zero. The superscript $*$ indicates dual space. The sign $\subset$
stands for proper inclusion of sets, and $\Z_{+}$ (respectively, $\Z_-$) denotes the set of positive 
(resp., negative) integer numbers.

 Let $\gk = \gk_0 \oplus \gk_1$ be a finite-dimensional Lie superalgebra,
  and $U(\gk)$ be the universal enveloping algebra of $\gk$. Fix
  a { Cartan subsuperalgebra} $\gh_\gk=(\gh_{\gk})_0 \oplus (\gh_{\gk})_1$ of $\gk$,
  i.e. a self-normalizing nilpotent subsuperalgebra. Then $(\gh_{\gk})_0$
   is a Cartan subsuperalgebra of $\gk_0$ and $(\gh_\gk)_1$ is the
   maximal subspace of $\gk_1$ on which $(\gh_\gk)_0$ acts nilpotently
    (see \cite {Sch}, \cite{PS2}).
By  $\Delta_{\gk} =(\Delta_{\gk})_0 \cup (\Delta_{\gk})_1  \subset (\gh_{\gk})^*_0 $ we denote the
     set of roots of $\gk$, and by $Q_{\gk}:= \Z \Delta_{\gk}$ we denote the root lattice of $\gk$. In this paper we consider a fixed Lie superalgebra $\gg$ for
which $\gh_\gg = (\gh_\gg)_0$ and $\gg_0 = \gg_0' \cplus \gg_0''$,
where    $\gg_0'$ is a reductive Lie algebra containing $\gh_\gg$
and $\gg_0''$ is a nilpotent Lie
    algebra. Set $\gh := \gh_\gg$, $\Delta := \Delta_\gg$, and $Q := Q_\gg$, and let $W$ denote the Weyl group of
    $\gg_0'$.

      A $\gg$-module $M$ is a
     {\it weight module}  if $M=\oplus_{\lambda\in \gh_0^*}
     M^{\lambda}$, where $M^{\lambda}:=\{ m \in M \; | \; h \cdot m=\lambda (h) m \}$.
     The space $M^{\lambda}$ is the {\it weight space of weight $\lambda$},
      and $\dim M^{\lambda}$ is the {\it multiplicity} of $M^{\lambda}$.
The {\it support} of $M$ is the set $\supp M := \{\lambda \in
\gh_0^* \; | \; M^{\lambda} \neq 0 \}$. The weight space
$U(\gg)^0$ of weight $0$ equals the supercentralizer of $\gh_0$ in
$U$.  $U:=U(\gg)$ is a weight $\gg$-module with respect to the
adjoint action of $\gg$. In the rest of the paper the term {\it
$\gk$-module} will be used as an abbreviation for a weight
$\gk$-module for any Lie superalgebra $\gk$.

Fix  a parabolic subsuperalgebra $\gp$ of $\gg$ which contains
$\gh$ and assume that $\gp$ admits a vector space decomposition
$\gp = \gu \oplus \ga$, where $\gu$ is the  nilradical of $\gp$,
and $\ga$ is a reductive Lie subalgebra of $\gp \cap \gg_0'$ with $\gh \subseteq \ga$.
Denote by $\gu^-$ the  nilradical of the parabolic subalgebra opposite to $\gp$. Then 
 $\gg = \gu^- \oplus \ga \oplus \gu$. 
If $S$ is an $\ga$-module put $M_{\gp}(S):=\Ind_{\gp}^{\gg}S$,
where $S$  is considered as a $\gp$-module with trivial action of
$\gu$. Suppose that the support of $S$ lies in a single $K \otimes
Q_{\ga}$-coset $u = \lambda + K \otimes Q_{\ga}$. Let $Z_{\gp}(S)$
be the sum of all $\gg$-submodules of $M_{\gp}(S)$ having trivial
intersection with $S$ ($M$ is the maximal $\gg$-module
$Z_{\gp}(S)$ included in $\oplus_{\mu \notin u} M_{\gp}(S)^{\mu}$,
see  Lemma 2.3 in \cite{DMP}). Set $L_{\gp}(S):= M_{\gp}(S) /
Z_{\gp}(S)$.  Corollary 2.4 in \cite{DMP} implies that $S$ is
simple if and only if $L_{\gp}(S)$ is simple.

 Finally, if $\gk \subseteq \gl$ are  Lie superalgebras,
and $P$ is an $\gl$-module,  denote by $H^0(\gk, P)$ the space of
$\gk$-invariant vectors in $P$.

\subsection{Definitions and properties of the twisted localization} For
the rest of the paper $\ga$ will be a  reductive Lie algebra, and
$\gd$ will be a Lie  superalgebra,  such that $ \ga \subseteq \gd
\subseteq \gg$. The results  obtained in this section will be
applied  to the cases  $\gd = \ga$, $\gd = \gg_0$, and $\gd =
\gg$. For every $\alpha \in \Delta_{\ga}$, let $h_{\alpha}$ be the
corresponding coroot, let $e_{\alpha}$ be a fixed generator of
$\ga^{\alpha}$, and let $f_{\alpha} \in \ga^{-\alpha}$ be defined
by the equality $[e_{\alpha}, f_{\alpha}] = h_{\alpha}$. Let
$\Gamma_{\ga} \subset \Delta_{\ga}$ be a commuting set of roots
which is a basis of $Q_{\ga}$, and let $F_{\Gamma_{\ga}}$ be the
multiplicative subset of $U$ generated by $\{f_{\alpha} \; | \;
\alpha \in \Gamma_{\ga} \}$. Since $F_{\Gamma_{\ga}}$ satisfies
Ore's localizability conditions, the localization
$U_{\Gamma_{\ga}}$ of $U$ with respect to $F_{\Gamma_{\ga}}$  is
well-defined. For details about this localization we refer the
reader to Section 4 in \cite{M}.

A $\gd$-module  $M$ will be called {\it $\Gamma_{\ga}$-injective}
(respectively, {\it $\Gamma_{\ga}$-bijective}) if $f_{\alpha}$
acts injectively (resp., bijectively) on $M$ for each $\alpha \in
\Gamma_{\ga}$. If $M$ is a simple, the condition that $f_{\alpha}$
acts injectively (resp., bijectively) on $M$ is equivalent to the
condition that the operator $f_{\alpha}: M^{\lambda} \to M^{\lambda - \alpha}$ is injective
(resp., bijective) for some $\lambda \in \supp M$, see Corollary 3.4
in \cite{DMP}. Let $M$ be a $\Gamma_{\ga}$-injective weight
$\gd$-module.  The $F_{\Gamma_{\ga}}${\it-localization} of $M$ is
$M_{F_{\Gamma_{\ga}}}:=U_{F_{\Gamma_{\ga}}}\otimes_{U}M$, where
$U_{F_{\Gamma_{\ga}}}$ is the localization of $U$ with respect to
$F_{\Gamma_{\ga}}$. There is a $\gd$-module injection $i : M \to
M_{F_{\Gamma_{\ga}}}$ which allows us  to identify  every
$\gd$-submodule $N$ of $M$ with  its image $i(N)$. It is clear
that $N_{F_{\Gamma_{\ga}}} = N$  if and only  if $N$ is
$\Gamma_{\ga}$-bijective.

If $\Gamma_{\ga} = \{ \alpha_1,...,\alpha_l\}$, the $l${\it
-parametric set of twisted automorphisms} of
$U_{F_{\Gamma_{\ga}}}$, $\Theta:K^l \times U_{F_{\Gamma_{\ga}}}
\mapsto  U_{F_{\Gamma_{\ga}}}$, is defined by the formula $$
\Theta_{(x_{1},\dots,x_{l})}(u):= \sum\limits_{0\leq
i_{1},\dots,i_{l}\leq N(u)} (^{x_{1}}_{i_{1}})\dots
(^{x_{l}}_{i_{l}})\,ad(f_{\alpha_1})^{i_{1}}\dots
ad(f_{\alpha_l})^{i_{l}}(u) \,f_{\alpha_1}^{-i_{1}}\dots
f_{\alpha_l}^{-i_{l}}, $$ where $(^{x}_{i}) :=
x(x-1)...(x-i+1)/i!$ for $x \in K$ and $i \in \Z_+ \cup \{ 0\} $.
The sum is well-defined, since, for every $u \in
U_{F_{\Gamma_{\ga}}}$, there exists a nonnegative integer $N(u)$
such that $ad(f_{\alpha_j})^{N(u)+1}(u)=0$ for  $j=1,2,..,l$. Note
that in the particular case when
 $(x_1,x_2,...,x_l) \in \Z^l$ we have $\Theta_{(x_{1},\dots,x_{l})}(u)= f_{\alpha_1}^{x_{1}}\dots f_{\alpha_l}^{x_{l}}\,u\,f_{\alpha_l}^{-x_{l}}\dots f_{\alpha_1}^{-x_{1}}$.

For $\mu = -(x_1 \alpha_1 + x_2 \alpha_2 +...+ x_l \alpha_l) \in K
\otimes Q_{\ga}$ and a $U_{F_{\Gamma_{\ga}}}$-module $N$, let $f_{\Gamma_{\ga}}^{\mu} N
$ be the module $N$ twisted by the automorphism
$\Theta_{(x_{1},\dots,x_{l})}$, i.e. the vector space
$N$   with $U_{F_{\Gamma_{\ga}}}$-action
defined by $u \cdot f_{{\Gamma_{\ga}}}^{\mu} v :=
f_{{\Gamma_{\ga}}}^{\mu} ( \Theta_{(x_{1},\dots,x_{l})}(u)\cdot
v)$, where for $v \in N$, $f_{{\Gamma_{\ga}}}^{\mu}  v$  stands for the element
$v$ considered as an element of the twisted module
$f_{{\Gamma_{\ga}}}^{\mu} N$. In particular,
if $v$ is a weight vector of weight $\eta$, then
$f_{{\Gamma_{\ga}}}^{\mu}  v$ is a weight vector of weight $\eta +
\mu$. Set $\Psi_{{\Gamma_{\ga}}}^{\mu}M :=
f_{{\Gamma_{\ga}}}^{\mu} M_{F_{\Gamma_{\ga}}}$. Note that
$\Psi_{{\Gamma_{\ga}}}^{\mu_1}M \simeq
\Psi_{{\Gamma_{\ga}}}^{\mu_2}M$ whenever $\mu_1 - \mu_2 \in
Q_{\ga}$ and therefore the $\gd$-module
$\Psi_{{\Gamma_{\ga}}}^{\mu}M$ is well-defined up to an
isomorphism for $\mu \in K\otimes Q_{\ga} / Q_{\ga}$. For any $\Gamma_{\ga}$-bijective
$\gd$-submodule $L$  of $\Psi_{{\Gamma_{\ga}}}^{\mu}M$, set
$\Phi_{{\Gamma_{\ga}}}^{-\mu}L:=f_{{\Gamma_{\ga}}}^{-\mu}L \cap M$
(the intersection is defined in $M_{F_{\Gamma_{\ga}}}$).

The map $\Psi_{{\Gamma_{\ga}}}^{\mu}$ from the lattice of
$\gd$-submodules of $M$ to the lattice of $\gd$-submodules of
$\Psi_{{\Gamma_{\ga}}}^{\mu}M$  plays central  role in the paper and this
section is devoted to its properties. The injectivity and
surjectivity of this map for the case when $\gg$ is a Lie algebra,
$\gd = \ga$ and $M = M_{\gp}(S)$ are also studied in \cite{BFL}. We
start with the observation that the image of
$\Psi_{{\Gamma_{\ga}}}^{\mu}$ always contains the set of all
$\Gamma_{\ga}$-bijective submodules of
$\Psi_{{\Gamma_{\ga}}}^{\mu}M$.

\begin{lemma} \label{Phi}
With the notations above, let $\mu \in K \otimes Q_{\ga}$ and  $N$ be a $\gd$-submodule of $\Psi_{{\Gamma_{\ga}}}^{\mu}M$ for which $N \simeq N_{F_{\Gamma_{\ga}}}$. Then $\Psi_{{\Gamma_{\ga}}}^{\mu}(\Phi_{{\Gamma_{\ga}}}^{-\mu}N) =N$.
\end{lemma}
\noindent {\bf Proof.} The statement is equivalent to $(M \cap
f_{{\Gamma_{\ga}}}^{-\mu}N)_{F_{\Gamma_{\ga}}} =
f_{{\Gamma_{\ga}}}^{-\mu}N$. Since $N$ is
$\Gamma_{\ga}$-bijective, so is $f_{\Gamma_{\ga}}^{-\mu}N$.
Therefore we have $(M \cap
f_{\Gamma_{\ga}}^{-\mu}N)_{F_{\Gamma_{\ga}}} \subseteq
(f_{\Gamma_{\ga}}^{-\mu}N)_{F_{\Gamma_{\ga}}} =
f_{\Gamma_{\ga}}^{-\mu}N$. To prove the inverse inclusion, observe
that for any $x \in f_{{\Gamma_{\ga}}}^{-\mu}N \subseteq
M_{F_{\Gamma_{\ga}}}$ there is  $f \in F_{\Gamma_{\ga}}$ for which
$f \cdot x \in M \cap f_{\Gamma_{\ga}}^{-\mu}N$. \hfill $\square$

\begin{lemma} \label{m1-m2}
 Let $M_1 \subset M_2$ be ${\Gamma_{\ga}}$-injective $\gd$-modules.

(i) The $\gd$-module $M_2 / M_1$ is ${\Gamma_{\ga}}$-injective if and only if  $M_1 = M_2 \cap (M_1)_{F_{\Gamma_{\ga}}}$. In particular, if $M_2 / M_1$ is simple, then it is $\Gamma_{\ga}$-injective if and only if  $(M_1)_{F_{\Gamma_{\ga}}} \neq (M_2)_{F_{\Gamma_{\ga}}}$.

(ii) If $M_2 / M_1$ is ${\Gamma_{\ga}}$-injective, then $(M_2)_{F_{\Gamma_{\ga}}} / (M_1)_{F_{\Gamma_{\ga}}}  \simeq (M_2/M_1)_{F_{\Gamma_{\ga}}}$. In particular $\Psi_{{\Gamma_{\ga}}}^\mu (M_2)/\Psi_{{\Gamma_{\ga}}}^\mu (M_1) \simeq \Psi_{{\Gamma_{\ga}}}^\mu (M_2/M_1)$  for every $\mu \in K \otimes Q_{\ga}$.
\end{lemma}
\noindent
{\bf Proof.} (i) The $\Gamma_{\ga}$-localizability of $M_2 /M_1$ fails if and only if $f_{\alpha}\cdot m_2 \in M_1$ for some $m_2 \in M_2 \setminus M_1$ and $\alpha \in \Gamma_{\ga}$. The latter is  equivalent to  $m_2 \in ((M_1)_{F_{\Gamma_{\ga}}} \cap M_2) \setminus M_1$, which implies the first part of the statement. If $M_2 /M_1$ is simple then the inclusions $M_1 \subseteq M_2 \cap (M_1)_{F_{\Gamma_{\ga}}} \subset M_2$ imply the second part. (ii) is a corollary of (i).
\hfill $\square$

\begin{lemma} \label{switch}
Let $M_1 \subset M_2 \subset M_3$ be a triple of
$\Gamma_{\ga}$-injective $\gd$-modules such that $M_2 / M_1$ is
$\Gamma_{\ga}$-injective and $(M_2)_{F_{\Gamma_{\ga}}} =
(M_3)_{F_{\Gamma_{\ga}}}$ (in particular $M_3 / M_2$ is not
$\Gamma_{\ga}$-injective). Set  $M_2':= M_3 \cap
(M_1)_{F_{\Gamma_{\ga}}}$. Then $M_1 \subset M_2' \subset M_3$ is
a triple of $\Gamma_{\ga}$-injective $\gd$-modules such that $M_3
/ M_2'$ is $\Gamma_{\ga}$-injective and $(M_1)_{F_{\Gamma_{\ga}}}
= (M_2')_{F_{\Gamma_{\ga}}}$.
\end{lemma}
\noindent {\bf Proof.} The inclusions $M_1 \subseteq M_2' \subset
(M_1)_{F_{\Gamma_{\ga}}}$ imply $(M_1)_{F_{\Gamma_{\ga}}} =
(M_2')_{F_{\Gamma_{\ga}}}$. The injectivity of $M_3 / M_2'$
follows from Lemma \ref{m1-m2}.

\hfill $\square$

The following lemma shows that the  the map
$\Psi_{{\Gamma_{\ga}}}^{\mu}$ (and in particular the
$F_{{\Gamma_{\ga}}}$-localization) commutes with the functors
$H^0$, $M_{\gp}$, and $L_{\gp}$.

\begin{lemma} \label{commut}
Let $\Gamma_{\ga} \subset \Delta_{\ga}$ be a set of commuting roots which is a basis of $Q_{\ga}$ and $R$ be an $\ga$-module whose support lies in a single $K \otimes Q_{\ga}$-coset.

(i) The following conditions are equivalent:

$\bullet$ the $\ga$-module $R$ is $\Gamma_{\ga}$-injective,

$\bullet$ the $\gg$-module $M_{\gp}(R)$ is $\Gamma_{\ga}$-injective,

$\bullet$ the $\gg$-module $L_{\gp}(R)$ is $\Gamma_{\ga}$-injective.

(ii) Let $R$ be $\Gamma_{\ga}$-injective. Then:

$\bullet$ $H^0(\gu,\Psi_{{\Gamma_{\ga}}}^{\mu} (L_{\gp}(R)))  \simeq  \Psi_{\Gamma_{\ga}}^{\mu} (H^0(\gu,  L_{\gp}(R))) \simeq \Psi_{\Gamma_{\ga}}^{\mu} R$,

$\bullet$ $L_{\gp}(\Psi_{{\Gamma_{\ga}}}^\mu (R)) \simeq \Psi_{{\Gamma_{\ga}}}^\mu(L_{\gp}(R))$,

$\bullet$ $M_{\gp}(\Psi_{{\Gamma_{\ga}}}^\mu (R)) \simeq \Psi_{{\Gamma_{\ga}}}^\mu(M_{\gp}(R))$,

\noindent
for all $\mu \in K \otimes Q_{\ga}$. In particular, $\Psi_{{\Gamma_{\ga}}}^\mu R$ is simple if and only if $\Psi_{{\Gamma_{\ga}}}^\mu(L_{\gp}(R))$ is simple.
\end{lemma}
\noindent {\bf Proof.} (i) If $R$ is $\Gamma_{\ga}$-injective then
$L_{\gp}(R)$ is also $\Gamma_{\ga}$-injective by Corollary 3.4 in
\cite{DMP}. A submodule of a $\Gamma_{\ga}$-injective module is
$\Gamma_{\ga}$-injective, and hence it remains to show  that
$M_{\gp}(R)$ and $L_{\gp}(R)$ are $\Gamma_{\ga}$-injective
whenever $R$ is $\Gamma_{\ga}$-injective. Suppose that $f_{\alpha}
\cdot m = 0$ for some  ${\alpha} \in \Gamma_{\ga}$ and $m \in
M_{\gp}(R)$. We may assume that $m \in M_{\gp}(R)^{\lambda} $ for
some $\lambda \in \gh^*$. By the PBW Theorem we have  $M_{\gp}(R)
\simeq U(\gu^-) \otimes S$, and in particular $m = \sum_{i = 1}^N
u_i \otimes  m_i$ for some $u_i \in  U(\gu^-)^{\lambda -
\lambda_i}$, $m_i \in S^{\lambda_i}$, and $\lambda_i \in \gh^*$.
If $u_i':= f_{\alpha}u_i - u_i f_{\alpha} \in U(\gu^-)^{\lambda -
\lambda_i + \alpha}$, then $\sum_{i=1}^N (u_i \otimes (f_{\alpha}
\cdot m_i) + u_i' \otimes m_i) = 0$. Let $i_0$ be such that
$\lambda_{i_0} + \alpha \neq \lambda_i$ for every $i = 1,...,N$.
Since $R$ is $\Gamma_{\ga}$-injective, $u_{i_0} \otimes
(f_{\alpha} \cdot m_{i_0}) \neq 0$, and thus $u_{i_0} \otimes
(f_{\alpha} \cdot m_{i_0}) + u''_{j_0} \otimes m_{j_0} =0 $, for
some $j_0$ and $u_{j_0}'' \in U(\gu^-)^{\lambda - \lambda_{j_0} +
\alpha}$. This implies that $\lambda_{i_0} + \alpha =
\lambda_{j_0}$ which contradicts to the choice of $i_0$. Hence
$M_{\gp}(R)$ is  $\Gamma_{\ga}$-injective.

We show next that $L_{\gp}(R)$ is $\Gamma_{\ga}$-injective. Since
$Z_{\gp}(R) \subseteq \oplus_{\mu \notin u} M_{\gp}(R)^{\mu}$,
where $u$ is the $K \otimes Q_{\ga}$-coset for which  $\supp R
\subset u$, we have  $(Z_{\gp}(R))_{\Gamma_{\ga}} \cap R = 0$.
Therefore $Z_{\gp}(R) \subseteq (Z_{\gp}(R))_{\Gamma_{\ga}} \cap
M_{\gp}(R) \subset M_{\gp}(R)$, and by the maximality of
$Z_{\gp}(R)$ we find that $Z_{\gp}(R) =(Z_{\gp}(R))_{\Gamma_{\ga}}
\cap M_{\gp}(R)$. Now Lemma \ref{m1-m2} implies that $L_{\gp}(R)$
is $\Gamma_{\ga}$-injective.

 (ii)  The proof of the statement is a modification of the proof of Lemma 13.2 in \cite{M}.  Set for simplicity $\gu^+ := \gu$. The inclusions  $[f_{\alpha}, \gu^{\pm}] \subseteq  \gu^{\pm}$ for every $\alpha \in {\Gamma_{\ga}}$ easily imply
\begin{equation} \label{upm}
f_{{\Gamma_{\ga}}}^{\mu} (\gu^{\pm} \cdot R_{F_{\Gamma_{\ga}}}) =  \gu^{\pm} \cdot (f_{{\Gamma_{\ga}}}^{\mu}  R_{F_{\Gamma_{\ga}}} )
\end{equation}
for all $\mu \in K\otimes Q_{\ga}$. Applying (\ref{upm}) for $\mu =
- \sum_{i=1}^k x_i \alpha_i$ and  positive integers  $x_i$ we
obtain $(U(\gu^{\pm}) \cdot R)_{F_{{\Gamma_{\ga}}}} = U(\gu^{\pm})
\cdot R_{F_{{\Gamma_{\ga}}}}$.  Therefore, using the PBW theorem
we have: $M_{\gp}(\Psi_{{\Gamma_{\ga}}}^\mu (R)) \simeq U(\gu^-)
\cdot (f_{{\Gamma_{\ga}}}^{\mu} R_{F_{{\Gamma_{\ga}}}}) =
f_{{\Gamma_{\ga}}}^{\mu} (U(\gu^-) \cdot R_{F_{{\Gamma_{\ga}}}}) =
f_{{\Gamma_{\ga}}}^{\mu}(U(\gu^-)  \cdot R)_{F_{{\Gamma_{\ga}}}}
=\Psi_{{\Gamma_{\ga}}}^\mu(M_{\gp}(R)) $ for all $\mu \in K
\otimes Q_{\ga}$. Finally, a double application of (\ref{upm}) yields  $$H^0(\gu, f_{{\Gamma_{\ga}}}^{\mu}
(L_{\gp}(R))_{F_{{\Gamma_{\ga}}}}) \simeq f_{{\Gamma_{\ga}}}^{\mu}
H^0(\gu, (L_{\gp}(R))_{F_{{\Gamma_{\ga}}}}) \simeq
f_{{\Gamma_{\ga}}}^{\mu} H^0(\gu,
L_{\gp}(R))_{F_{{\Gamma_{\ga}}}},$$ which completes the proof of
all statements in (ii). \hfill $\square$

The following theorem is our main result in this section.

\begin{theorem} \label{Psi}
Let $M$ be a  $\Gamma_{\ga}$-injective $\gd$-module of finite length.

(i) There exists a composition series
\begin{equation} \label{JH2}
0 = M_1^0 \subset M_1^1 \subset ... \subset M_{r_1}^1 \subset M_1^2 \subset ... \subset M_{r_2}^2 \subset ... \subset M_1^t \subset ... \subset M_{r_t}^t = M
\end{equation}
of $M$ such that

$\bullet$ $(M_{i_1}^{j_1})_{F_{\Gamma_{\ga}}} =(M_{i_2}^{j_2})_{F_{\Gamma_{\ga}}} $ if and only if $j_1 = j_2$;

$\bullet$ all modules $M_j:= M_{r_j}^j/M_{r_{j-1}}^{j-1}$ are $\Gamma_{\ga}$-injective.

(ii) The modules $M_j$ are indecomposable.

(iii) For every $\mu \in K \otimes Q_{\ga}$ the $\gd$-module $\Psi_{\Sigma}^\mu (M)$ admits a filtration $0 \subset \Psi_{\Sigma}^\mu (M_1^1) = ... = \Psi_{\Sigma}^\mu (M_{r_1}^1) \subset \Psi_{\Sigma}^\mu(M_1^2) = ... = \Psi_{\Sigma}^\mu(M_{r_2}^2) \subset ... \subset \Psi_{\Sigma}^\mu(M_1^t) = ...= \Psi_{\Sigma}^\mu(M_{r_t}^t)  = \Psi_{\Sigma}^\mu (M)$ with successive nontrivial quotients $\Psi_{\Sigma}^\mu(M_j)$.

\end{theorem}
\noindent {\bf Proof.} (i) Starting with an arbitrary composition
series of $M$ we enumerate the successive submodules so that the
first condition is satisfied. In order to obtain
$\Gamma_{\ga}$-injective quotients $M_j$ we  modify the series by
consequently applying Lemma \ref{switch} for $j = t, t-1,...,2$. As
a result we obtain a series for which $M_j$ are indecomposable for
$j>1$. Finally, $M_1 \simeq M_{r_1}^1$ is always
$\Gamma_{\ga}$-injective as a submodule of $M$ and thus the
assertion is proved.

(ii) Lemma \ref{m1-m2} implies that among the simple quotients
$M_1^i/M_1^{i-1}$ only $M_1^1 \simeq M_1^i/M_1^{0}$ is
$\Gamma_{\ga}$-injective. Every simple submodule of $M_1 \simeq
M_{r_1}^1$ is $\Gamma_{\ga}$-injective and therefore is isomorphic
to $M_1^1$. In particular, every submodule $X$ of $M_1$ has a
simple submodule isomorphic to $M_1^1$. Thus
$(M_1^1)_{F_{\Gamma_\ga}} \simeq  X_{F_{\Gamma_{\ga}}} \simeq
(M_1)_{F_{\Gamma_{\ga}}}$. Suppose now $M_1 \simeq  A \oplus B$
for some $\gd$-modules $A$ and $B$. Then $A_{F_{\Gamma_\ga}}
\simeq B_{F_{\Gamma_\ga}} \simeq (M_1)_{F_{\Gamma_{\ga}}}$ on one
hand, and $(M_j)_{F_{\Gamma_{\ga}}} \simeq  A_{F_{\Gamma_{\ga}}}
\oplus B_{F_{\Gamma_{\ga}}}$ on the other, which is impossible.
Therefore $M_1$ is indecomposable. In order to show that $M_j$ is
indecomposable for $j>1$ we apply induction on $j$ observing that
all submodules of the $\Gamma_{\ga}$-injective module $M_j:=
M_{r_j}^j/M_{r_{j-1}}^{j-1}$ have the same
$\Gamma_{\ga}$-localization.

(iii)  follows directly from Lemma  \ref{m1-m2}, (ii).
\hfill $\square$

We call a series (\ref{JH2}) with the properties of Theorem \ref{Psi}, (i), a  {\it $\Gamma_{\ga}-injective$ composition series} of $M$.

\section{Coherent Families of Weight Modules}

In the rest of the paper $\gg$ will be one of the following  Lie
superalgebras: $\gg \gl(m|n)$, $\gs \gl (m|n)$, for $m \neq n$, $\gp \gs \gl (m|m)$, $\go \gs \gp (2|n)$ (where $n = 2q$), $\pp (m)$, $\gs {\pp} (m)$,  
or {\bf W}$(n)$. We will also assume  that all weight modules considered have
finite weight multiplicities. In this section we define a $(\gg,\ga)$-coherent
family and show that every simple $\gg$-module $M$ is a submodule
of a unique semisimple coherent family ${\cal E}_{\ga}(M)$.

We start by recalling some basic facts about the root structure
of the listed above Lie superalgebras. More details can be found in \cite{P}. For all  Lie superalgebras $\gg$  we have $\gh = \gh_0$
and $\gh_1=0$. Let $\Delta = \Delta_0 \cup \Delta_1$ be the set of roots of
$\gg$ and let $\Delta_1 =\Delta_1^+ \cup \Delta_1^- $, where $\Delta_0$ and $\Delta_1^{\pm}$ are written explicitly below.

$\bullet$ If $\gg = \gg \gl(m|n)$, $\gg = \gs \gl (m|n)$ for $m \neq n$, or $\gg = \gp \gs \gl (m|m)$ for $m=n$ then 

\noindent
$\gg_0 = \gg_0' = \gg \gl (m) \oplus \gg \gl (n)$ for  $\gg = \gg \gl(m|n)$, $\gg_0 = \gg_0' \simeq \gs \gl (m) \oplus \gs \gl (n) \oplus K$ for  $\gg = \gs \gl(m|n)$, and $\gg_0 = \gg_0' \simeq \gs \gl (m) \oplus \gs \gl (m)$ for  $\gg = \gp \gs \gl(m|m)$. The vector space $\gh_{\gg \gl (m|n)}^*$ has a basis $\{ \varepsilon_1,...,  \varepsilon_m, \delta_1,...,\delta_n\}$ and $\gh_{\gs \gl (m|n)}^* \simeq \gh_{\gg \gl (m|n)}^* / (K(\sum_i \varepsilon_i - \sum_j \delta_j))$. For the images of $\varepsilon_i$ and $\delta_j$ in $\gh_{\gs \gl (m|n)}^*$ we use the same letters. Furthermore, $\gh_{\gp \gs \gl (m|m)}^* \simeq \{ \sum_i k_i \varepsilon_i + \sum_i l_i \delta_i \in \gh_{\gs \gl (m|m)}^* \; | \; \sum_i k_i + \sum_i l_i = 0\} \subset \gh_{\gs \gl (m|m)}^*$. In all cases $\Delta_0 = \{ \varepsilon_i - \varepsilon_j, \delta_k - \delta_l \; | \; 1\leq i \neq j \leq m, 1 \leq k \neq l \leq n \} $, 
$\Delta_1^{\pm} = \{ \pm (\varepsilon_i - \delta_l) \; | \; 1\leq i \leq m, 1 \leq l \leq n \} $.

$\bullet$ If $\gg = \go \gs \gp (2|n)$, $n = 2q$ then 

\noindent
$\gg_0 = \gg_0' = \go (2) \oplus \gs \gp (n) \simeq K \oplus \gs \gp (n)$. The vector space  $\gh^*$ has a basis $\{ \varepsilon_1, \delta_1,...,\delta_q\} $,
$\Delta_0 = \{ \pm 2 \delta_k, \pm \delta_k \pm \delta_l\; | \; 1 \leq k \neq l \leq q \} $, 
$\Delta_1^{\pm} = \{ \pm (\varepsilon_1 + \delta_k),\pm (\varepsilon_1 - \delta_k)  \; | \; 1 \leq k \leq q \} $.

$\bullet$ If $\gg = \pp (m)$, or $\gg = \gs {\pp} (m)$, then 

\noindent
$\gg_0 = \gg_0' = \gg \gl (m)$ for  $\gg =  \pp (m)$ and $\gg_0 = \gg_0' = \gs \gl (m)$ for  $\gg = \gs {\pp} (m)$.
The vector space  $\gh_{p(m)}^*$ has a basis $\{ \varepsilon_1, ..., \varepsilon_m\} $ and  $\gh_{\gs p(m)}^* \simeq \gh_{p(m)}^* / (K(\sum_i \varepsilon_i))$.  For the images of $\varepsilon_i$ in $\gh_{\gs p (m)}^*$ we use the same letters.
$\Delta_0 = \{ \varepsilon_i - \varepsilon_j \; | \; 1 \leq i \neq j \leq m \} $, 
$\Delta_1^{+} = \{ \varepsilon_i + \varepsilon_j, 2\varepsilon_i  \; | \; 1 \leq i \neq j \leq m \} $, 
$\Delta_1^{-} = \{ -(\varepsilon_i + \varepsilon_j)  \; | \; 1 \leq i \neq j \leq m \} $.

$\bullet$ If $\gg =$ {\bf W}$(n)$ then

\noindent
$\gg_0' = \gg \gl (n)$ and $\gg_0''$ is a nilpotent Lie algebra of dimension $n 2^{n-1} - n^2$. The vector space  $\gh^*$ has a basis $\{ \varepsilon_1, ..., \varepsilon_n\} $. 
$\Delta_0 = \{ \varepsilon_{i_1} +...+ \varepsilon_{i_k}, \varepsilon_{i_1} +...+ \varepsilon_{i_l} - \varepsilon_{j} \; | \; k = 2t, l = 2u+1\} $, 
$\Delta_1 = \{ \varepsilon_{i_1} +...+ \varepsilon_{i_k}, \varepsilon_{i_1} +...+ \varepsilon_{i_l} - \varepsilon_{j} \; | \; k = 2t+1, l = 2u\} $,

All Lie superalgebras admit a natural $\Z$-grading. If $\gg \neq ${\bf W}$(n)$ then  $\gg = \gg^{-1} \oplus \gg^0 \oplus \gg^{1}$, where $\gg^{\pm 1}:= \oplus_{\alpha \in \Delta_1^{\pm}} \gg^{\alpha}$ and $\gg^0 := \gg_0$. For $\gg =$ {\bf W}$(n)$ we consider the $\Z$-grading {\bf W}$(n) :=  \oplus_{k=-1}^{n-1}$ {\bf W}$^k$, where {\bf W}$^k := \oplus_{\sum t_i = k}${\bf W}$(n)^{\sum t_i \varepsilon_i}$. We also set  $W_{\geq k} := \oplus_{i \geq k} ${\bf W}$^i$.

If $\gc$ is a reductive Lie algebra,
a $\gc$-module $R$ is called {\it torsion-free} if for every
$\alpha \in \Delta_{\gc}$, every  $x \in \gc^{\alpha}$ acts
bijectively on $R$. The following description of all simple weight $\gg$-modules follows from a set of results in \cite{DMP} and is crucial in our paper.

\begin{proposition} \label{dmp}
Let $M$ be a simple $\gg$-module. There exists a parabolic subsuperalgebra $\gq \subset \gg$  with nilradical 
$\gn_{\gq}$ for which $\gq = \gc \oplus \gn_{\gq}$ and $\gc$ is a reductive Lie subalgebra of $\gg_0'$ containing $\gh$, such that 
$M \simeq L_{\gq}(R)$ for some torsion-free $\gc$-module $R$. The list of possible Lie algebras $\gc$ is provided below.

$\bullet$ If $\gg = \gg \gl(m|n)$ (respectively, $\gg = \gs \gl (m|n)$), $m \neq n$, then $\gc$ is isomorphic to $\gg \gl (m_1,...,m_k,n_1,...,n_l)$ (resp., $\gs \gl (m_1,...,m_k,n_1,...,n_l)$ ), for some positive integers $m_i$ and $n_j$ with $\sum_i m_i = m$ and $\sum_j n_j = n$.

$\bullet$ If $\gg = \gp \gs \gl(m|m)$ then $\gc$ is isomorphic to $\gs \gl (m_1,...,m_k) \oplus \gs \gl(m_1',...,m_l')$, for some positive integers $m_i$ and $m_j'$ with $\sum_i m_i = \sum_j m_j' = m$. 

$\bullet$ If $\gg = \go \gs \gp (2|n)$, $n = 2q$ then $\gc$ is isomorphic to $k \oplus \gg \gl (q_1',...,q_k') \oplus \gs \gp (2q'')$, for some positive integers $q_i'$ and $q''$ with $\sum_i q_i' = q'$ and $q' + q'' = q$.

$\bullet$ If  $\gg = \gs \pp (m)$ (respectively, $\gg = {\pp} (m)$) then $\gc$ is isomorphic to $\gs \gl (m_1,...,m_k)$ (resp., $\gg \gl (m_1,...,m_k)$), for some positive integers $m_i$ with $\sum_i m_i = m$.

$\bullet$ If  $\gg = $ {\bf W}$(n)$ then $\gc$ is isomorphic to $\gg \gl (n_1,...,n_k)$,  for some positive integers $n_i$ with $\sum_i n_i = n$.
 
(For a $k$-tuple of positive integers $(m_1,...,m_k)$ with $\sum_i m_i = m$ by  $\gs \gl(m_1,...,m_k)$ we denote  $\{ (g_1 +...+ g_k) \in \gg \gl(m_1,...,m_k):=\gg \gl (m_1) \oplus ... \oplus \gg \gl (m_k)\; | \; \sum_i \tr(g_i) = 0\}$.)
\end{proposition}
\noindent {\bf Proof.} Assume for simplicity that $\gg$ is simple (for $\gg = \gg \gl (m|n)$ and $\gg = \pp(m)$ the statement easily follows from the cases of $\gg = \gs \gl (m|n)$ and $\gg = \gs \pp (m)$ respectively).  Following \cite{DMP}, we recall some definitions which will be used in this proof only. For a finite dimensional Lie superalgebra $\gk$ and a simple $\gk$-module $R$ denote by $inj(R)$ the set of all roots $\alpha \in (\Delta_{\gk})_0$ such that $x$ acts injectively on $R$  for some $x \in (\gk_0)^{\alpha}$. Let $C_{inj(R)}$ be the cone in $Q_{\gk}$ generated by $inj(R)$,  and let $C_R:= \{\lambda \in Q_{\gk} \; | \; m\lambda \in C_{inj(R)},$ for some $m \in \Z\}$ be the {\it saturation} of $C_{inj(R)}$. Theorem 3.6 and Corollary 3.7 in \cite{DMP} imply that $M$ is isomorphic to the unique simple quotient $L_{T}(\Omega)$ of $M_{\gg_T^0 \oplus \gg_T^+}(\Omega)$ where $\gg = \gg_T^- \oplus \gg_T^0 \oplus \gg_T^+$, $\gg_T^{\pm}:=\oplus_{l(\alpha) \in \Z_{\pm}}\gg^{\alpha}$, $\gg_T^{0}:=\oplus_{l(\alpha) = 0}\gg^{\alpha}$ for some linear map $l : Q \to \Z$ and a simple $\gg_T^0$-module $\Omega$ for which $C_{\Omega} = Q_{\gg_T^0}$. Furthermore, we use the classification of all possible Lie superalgebras $\gg_T^0$ provided in Section 7 of \cite{DMP}. We conclude that $(\gg_T^0)^1 = 0$ except for the following couples $(\gg, \gg_T^0)$: $(\gp \gs \gl (m|m), \gs \gl (m|m))$, $(\gs {\pp} (m),\gs {\pp} (m)) $ and $(\gs {\pp} (2m),  \gs \gl (m|m))) $. In the remaining three cases we observe that $\Omega \simeq L_{(\gg_T^0)^0 \oplus (\gg_T^0)^1}(S)$, where $S := H^0((\gg_T^0)^1, \Omega)$ is a torsion-free $(\gg_T^0)^0$-module. Therefore $M \simeq L_{(\gg_T^0)^0 \oplus((\gg_T^0)^1 + \gg_T^+)}(S)$ and we complete the proof. \hfill $\square$

\bigskip
\noindent
{\bf Remark.} The fact that every simple $\gg$-module is  isomorphic to a quotient of a module induced from a parabolical subalgebra,  whose reductive part is a Lie algebra, is the main reason to restrict our attention to the listed Lie superalgebras and to study the structure of the $\gg$-modules $L_{\gp}(S)$ and $M_{\gp}(S)$.

\bigskip

 Let $\ga \simeq \gs \oplus \gz \simeq \ga_1 \oplus ... \oplus \ga_k \oplus \gz $,
where $\ga_i$ are simple Lie algebras and $\gz$ is the center of
$\ga$. For every weight $\eta \in \gh^* \simeq (
\oplus_{i=1}^{k}\gh_{\ga_i}^* )\oplus \gh_{\gz}^*$ we write
$\eta^{\ga_i}$,  $\eta^{\gz}$, and $\eta^{\gs}$ for the
projections of $\eta$  on $\gh_{\ga_i}^*$, $\gh_{\gz}^*$, and
$\gh_{\gs}^*$ respectively, i.e. $\eta = \eta^{\gs} + \eta^{\gz} =
\sum_{i=1}^{k} \eta^{\ga_i} + \eta^{\gz}$.  Fix a basis $B$  of $\Delta$ such that $B \subset \Delta_{\gp}$, and set
$B_{\ga}:=B\cap\Delta_{\ga}$. For any $\lambda\in {\gh}^{*}$, let $M_{B_{\ga}}(\lambda)$ (respectively, $M_{B}(\lambda)$) be the Verma module
with $B_{\ga}$-highest (resp., $B$-highest) weight $\lambda$
and  $L_{B_{\ga}}(\lambda)$  (resp., $L_{B}(\lambda) \simeq  L_{\gp}(L_{B_{\ga}}(\lambda))$) be its unique simple quotient.  Set $T^{*}_{\ga}:= \gh^* /Q_{\ga}$,  $T^{*}_{\gs}:= (K \otimes Q_{\ga}) / Q_{\ga} \simeq \gh_{\gs}^* /Q_{\gs}$ and $T^{*} := \gh^*/Q_{\gg}$.  For $U \subset \gh^*$ and an $\gh$-module $M$ define $M[U]:= \oplus_{\eta \in U}M^{\eta}$. The elements of $T^{*}$ and $T^{*}_{\ga}$ will be considered as subsets of $\gh^*$. In particular, for a $\gg$-module  $M$ and $t \in T^*$,  $M[t]$ is a $\gg$-submodule of $M$.  Set ${\gh}^{*}_{rel}:={\gh}^{*}/K\otimes Q_{\ga} \simeq \gz$.
   For any ${\ga}$-module $M$ and $t \in {\gh}^{*}_{rel}$,
put  $\deg_{ \ga }\,M[t]:= {\sup}_{\mu\in t}\,{\dim}\, M^{\mu}$.
The function  $\deg_{\ga}:t\mapsto \deg_{\ga}\,M[t]$ is called the
{\it $\ga$-relative degree of $M$}.  We  will call
$M$  $\ga${\it -bounded} if $\deg_{\ga}\,M[t]<\infty$ for any
$t\in \gh_{rel}^{*}$ (In \cite{M}, Mathieu uses the term ``admissible'' weight module, but to avoid confusion with Harish-Chandra modules of finite type we prefer to use the term ``bounded'' weight module).

\begin{lemma} \label{fin-len}
Let $S$ be an $\ga$-bounded $\ga$-module, Then $M_{\gp}(S)$ (and thus also $L_{\gp}(S)$)  is an $\ga$-bounded  $\gg$-module of finite length.

\end{lemma}
\noindent {\bf Proof.} Being an $\ga$-bounded module, $S$ has
finite length (Lemma 3.3 in \cite{M}), and since $M_{\gp}$ is an
exact functor  we may  assume that $S$ is simple. The
admissibility of $M_{\gp}(S)$ follows from the the admissibility
of $S$ and the isomorphism $M_{\gp}(S) \simeq U(\gu^-) \otimes S$.
Indeed, if $t \in \gh_{rel}^*$ and $\Delta_{\gu^-}$ is the root
system of $\gu^-$, the set $$X(t):= \{ (\beta_1, ..., \beta_r)\; |
\; \sum_{i=1}^r \beta_i \in \supp (M_{\gp}(S)[t]) - \supp S,
\beta_i \in \Delta_{\gu^-}\}$$
 is finite. Then $N(t):= \sum_{(\beta_1,...\beta_r) \in X(T)}\dim (U(\gu^-)^{\sum
\beta_i}) < \infty$ and therefore\\ $\deg_{\ga}(M_{\gp}(S)[t]) \leq
N(t) \deg_{\ga}S$ for every $t \in \gh^*_{rel}$. In order to show that $M_{\gp}(S)$
has finite length recall that every finitely generated weight module with finite weight
 multiplicities over a reductive Lie algebra has finite length (Theorem 4.21 in \cite{F}). Since for any $s \in S$, $M_{\gp}(S) \simeq U(\gu^-) \otimes S \simeq U\cdot (1 \otimes s)$ is finitely generated $\gg_0'$-module, , $M_{\gp}(S)$ has finite length. \hfill $\square$

 Assume that $M$ is a $\gg$-module for which $\supp M$ is included in a single $Q$-coset. Then $M[\lambda + Q_{\ga}] = M[\lambda + K \otimes Q_{\ga}]$ for every
$\lambda \in \gh^*$, and therefore whenever we consider $M[t]$ we may assume that $t \in T_{\ga}^*$. In this case set $\supp_{T_{\ga}^*}M := \{ t \in T_{\ga}^* \, | \, M[t] \neq 0 \} \subset Q_{\gg}/Q_{\ga}$. The partial order of $Q_{\gg}$ determined by $B$ defines naturally a partial order on $ \supp_{T_{\ga}^*}M $. If  $M$ is ${\ga}$-bounded, its {\it
$\ga$-essential support}
is ${\supp}_{\ga-ess}\,M:=
\cup_{t\in T_{\ga}^{*}}
\,\{\lambda\in t\vert\,{\dim}\,M^{\lambda}=\deg_{\ga}\,M[t]\}$.
Then $M$ is called {\it strictly }$\ga${\it -bounded} (or  just {\it bounded} in the special case of $\ga = \gg_0'$)
if ${\supp}_{\ga-ess}\,M\cap t$ is Zariski dense in the
vector space $K\otimes t$
for any $t\in \supp_{T_{\ga}^*}M$.

The motivation of introducing the above terminology is contained in the following slight modification of Lemma 13.1 in \cite{M}. We skip the proof which essentially repeats the arguments of Mathieu's proof.

\begin{lemma} \label{stradm}
Let $S$ be a simple $\ga$-module.

(i) The following are equivalent:

$\bullet$ $S$ is strictly $\ga$-bounded $\ga$-module,

$\bullet$ $L_{\gp}(S)$ is strictly $\ga$-bounded  $\gg$-module,

$\bullet$ $M_{\gp}(S)$ is strictly $\ga$-bounded $\gg$-module.

(ii) If $S$ is $\ga$-bounded, then it is strictly bounded if
and only if there exists a basis $\Gamma_{\ga}$ of $\Delta_{\ga}$
which consist of commuting roots such that $S$ is
$\Gamma_{\ga}$-injective. \hfill $\square$

\end{lemma}

Let  $\lambda \in \gh^*$. Lemma \ref{commut} implies that if one of the modules $L_{B_{\ga}}(\lambda)$, $L_{B}(\lambda)$, or $M_{\gp}(\lambda)$ is  $\Gamma_{\ga}$-injective then so are the other two. In that case we say that $\lambda$ is {\it $\Gamma_{\ga}$-injective} (or  {\it $(B, \Gamma_{\ga})$-injective} when $B$ is not fixed). We say also that $\lambda= \sum_{i=1}^{k} \lambda^{\ga_i} + \lambda^{\gz}$ is $\ga${\it -partially finite} if  $\lambda^{\ga_i}$ is dominant integral for some $i$. Furthermore, we call a  $\gd$-module $L$ of finite length  $\ga${\it -partially finite} if every simple subquotient of $L$ is a highest weight module with an $\gd$-partially finite highest weight.

\begin{lemma} \label{inj-fin}
(i) Let $L_{B_{\ga}}(\lambda)$ be an bounded module. Then $\lambda$ is  $\Gamma_{\ga}$-injective if and only if $\lambda$ is not an $\ga$-partially finite weight.

(ii) Let $L_1 \subseteq L_2$ be  $\Gamma_{\ga}$-injective $\gd$-modules such that $L_2 / L_1$ has finite length. Then, if  $(L_1)_{F_{\Gamma_{\ga}}} = (L_2)_{F_{\Gamma_{\ga}}}$, the $\gd$-module $L_2 / L_1$ is  $\ga$-partially finite.
\end{lemma}
\noindent
{\bf Proof.} (i) It is clear that  no $\lambda^{\ga_i}$ is dominant integral for a $\Gamma_{\ga}$-injective weight $\lambda$. The inverse direction follows from Lemma \ref{stradm}, (ii). (ii) Since  $(L_1)_{F_{\Gamma_{\ga}}} = L_{F_{\Gamma_{\ga}}} = (L_2)_{F_{\Gamma_{\ga}}}$ for every $\gd$-module $L$ such that $L_1 \subseteq L \subseteq L_2$, it is enough to consider the case when $L_2/L_1$ is simple. Then the statement follows from part (i) and Lemma \ref{m1-m2}, (i).
\hfill $\square$

Following \cite{M}, if  $d$ is a positive integer, an $\gs$-{\it
coherent family} $\cal M$ {\it of degree} $d$ is an  $\gs$-module
satisfying the following two conditions:\\ (i) $\dim {\cal
M}^{\lambda} = d$ for any $\lambda \in \gh_{\gs}^*$;  \\ (ii) the
{trace} function $\lambda  \mapsto \tr(u)\vert_{{\cal
M}^{\lambda}}$ is polynomial in $\lambda \in \gh_{\gs}^*$ for any
fixed $u \in U(\gs)^0$.

 We define an $\ga$-module ${\cal M}$ to be an {\it $\ga$-coherent family} if ${\cal M}[h]$ is $\gs$-coherent family for every $h \in \gh^*_{rel}$. In particular, every $\ga$-coherent family is $\ga$-bounded and is a (possibly infinite) direct sum of $\gs$-coherent families which are $\ga$-modules. A {\it $(\gg, \ga)$-coherent family} is a  $\gg$-module ${\cal M}$ which is an $\ga$-coherent family. A $(\gg, \ga)$-coherent family ${\cal M}$ is said to be {\it irreducible} if ${\cal M}[t]$ is a simple $\gg $-module for some $t \in T^*$. ${\cal M}$ is  {\it semisimple} if
${\cal M}[t]$ is a semisimple $\gg$-module for any $t \in T^*$.

For an  $\ga$-bounded $\gg$-module $M$, a $(\gg,\ga)$-{\it
coherent extension} ${\cal E}$ of $M$ is a $(\gg,\ga)$-coherent
family which contains $M$ as a subquotient and $\deg_{\ga} M =
\deg_{\ga}{\cal E}$. For any $\Gamma_{\ga}$-injective $\ga$-module $R$, we set ${\cal E}_{\Gamma_{\ga}}(R) : = \oplus_{\mu \in T_{\gs}^*} \Psi_{\Gamma_{\ga}}^{\mu}(R) $.

\begin{proposition} \label{ext}

Let $S$ be a strictly $\ga$-bounded simple $\ga$-module, and
let  $\Gamma_{\ga} \subset \Delta_{\ga}$ be such that $S$ is
$\Gamma_{\ga}$-injective (the existence of $\Gamma_{\ga}$ follows
from Lemma \ref{stradm}). Then  ${\cal E}_{\Gamma_{\ga}}(S)$,
${\cal E}_{\Gamma_{\ga}}(L_{\gp}(S))$, and ${\cal
E}_{\Gamma_{\ga}}(M_{\gp}(S))$ are $(\gg,\ga)$-coherent extensions
of $S$, $L_{\gp}(S)$, and $M_{\gp}(S)$ respectively. Moreover,
$L_{\gp}({\cal E}_{\Gamma_{\ga}}(S)) \simeq {\cal
E}_{\Gamma_{\ga}}(L_{\gp}(S))$ and $M_{\gp}({\cal
E}_{\Gamma_{\ga}}(S)) \simeq {\cal E}_{\Gamma_{\ga}}(M_{\gp}(S))$.
\end{proposition}
\noindent {\bf Proof.} For every $\Gamma_{\ga}$-injective
$\ga$-module $R$ we have  $\deg_{\ga} R = \deg_{\ga}
\Psi_{\Gamma_{\ga}}^{\mu}(R) $ for any $\mu \in T_{\gs}^*$ which
implies the first part of the proposition. The second part is a
corollary of Lemma \ref{commut}, (ii). \hfill $\square$

\begin{proposition} \label{cohext}
Let $M$ be a simple $\gg$-module isomorphic to  $L_{\gp}(S)$ for some strictly $\ga$-bounded $\ga$-module $S$.

(i) There exists a unique semisimple $(\gg, \ga)$-coherent extension
${\cal E}_{\ga}(M)$ of $M$ isomorphic to  $(L_{\gp}({\cal E}_{\Gamma_{\ga}}(S)))^{ss} \simeq ({\cal E}_{\Gamma_{\ga}}(M))^{ss} $ for some $\Gamma_{\ga} \subset \Delta_{\ga}$ such that $S$ is $\Gamma_{\ga}$-injective.

(ii) The $(\gg, \ga)$-coherent family ${\cal E}_{\ga}(M)$ is irreducible.
For any simple strictly $\ga$-bounded $\gg$-submodule $N$ of ${\cal E}_{\ga}(M)$,
we have ${\cal E}_{\ga}(M)\simeq {\cal E}_{\ga}(N)$.

(iii) Let $S$ be a torsion-free $\ga$-module. For any basis $D$ of $\Delta$, such that $D \subset \Delta_{\gp}$, there is a strictly $\ga$-bounded highest weight $\gg$-submodule $L_D(\lambda)$  of ${\cal E}_{\ga}(M)$.  Moreover, there exists a set $\Sigma_{\ga} \subset \Delta_{\ga}$ of commuting roots which is a basis of $Q_{\ga}$ such that:

$\bullet$ $ \Psi_{\Sigma_{\ga}}^\mu (L_{D_{\ga}}(\lambda)) \simeq S$, where $D_{\ga} := D \cap \Delta_{\ga}$,

$\bullet$ $ \Psi_{\Sigma_{\ga}}^\mu (L_{D}(\lambda)) \simeq M$,

$\bullet$ $ \Psi_{\Sigma_{\ga}}^\mu (M_{\gp}(\lambda)) \simeq M_{\gp}(S)$,
where $\mu \in \gh^*$ is any weight for which $\mu + \lambda \in \supp S$.

\end{proposition}
\noindent
{\bf Proof.} (i) Proposition \ref{ext} implies that  ${\cal E}_{\Gamma_{\ga}}(M) \simeq L_{\gp}({\cal E}_{\Gamma_{\ga}}(S))$ is  a  coherent extension of $M$. Using Lemma \ref{fin-len} we see that $L_{\gp}({\cal E}_{\Gamma_{\ga}}(S))[t]$ has finite length for every $t \in T^*$ and thus we may define a semisimple coherent extension  ${\cal E}_{\ga}(M) := \oplus_{t \in T^*}(L_{\gp}({\cal E}_{\Gamma_{\ga}}(S))[t])^{ss}$ of $M$. Let ${\cal E}'$ be another semisimple coherent extension of $M$. The $\ga$-modules ${\cal E}'$ and ${\cal E}_{\ga}(M)$ have the same trace on the Zariski dense  subset ${\supp}_{\ga-ess} M \cap t$ of $t \otimes K$ for every $t \in \supp_{T^*_{\ga}}M$. Since a  semisimple $\ga$-module is determined uniquely by its trace (Lemma 2.3 in \cite{M}) we have that  ${\cal E}'[s]$ and ${\cal E}_{\ga}(M)[s]$ are isomorphic $\ga$-modules for every $s \in \gh^*_{rel}$. It remains to show that they are isomorphic $\gg$-modules. If $\lambda_0 \in \supp S$, then $t_0 = \lambda_0 + Q_{\ga}$ is a minimal element in $\supp_{T_{\ga}^*}M$ (in particular $M[t_0] \simeq S$).  Thus $s_0:= \lambda_0 + \mu + Q_{\ga}$ is a minimal element in $\supp_{T_{\ga}^*} {\cal E}'[s] = \supp_{T_{\ga}^*} {\cal E}_{\ga}(M)[s]$ for any $\mu \in K \otimes Q_{\ga}$ such that $\lambda_0 + \mu  \in s$. Since  ${\cal E}'[s_0] \simeq {\cal E}_{\ga}(M)[s_0]$ has finite length we can choose isomorphic simple $\ga$-submodules $L'$ and $N'$ of ${\cal E}'[s_0]$ and ${\cal E}_{\ga}(M)[s_0]$ respectively. Then $ {\cal E}'[s_0] \simeq L_{\gp}(L') \oplus L''$ and  $ {\cal E}_{\ga}(M)[s_0] \simeq L_{\gp}(N') \oplus N''$ for some semisimple $\gg$-modules $L''$ and $N''$ with isomorphic $\ga$-composition factors. We prove that $ {\cal E}'[s]$ and  $ {\cal E}_{\ga}(M)[s]$ are isomorphic $\gg$-modules by induction on the length of $ {\cal E}_{\ga}(M)[s]$ (which is finite by Lemma \ref{fin-len}).

(ii) Let $r_0, t_0 \in T^*_{\ga}$ be such that $M[t_0] = S$ and
$N[r_0] = H^0(\gu , N)$. Then $N[r_0]$ is a strictly
$\ga$-bounded module. By Lemma \ref{stradm} we find a basis
$\Gamma_{\ga}'$ of $\Delta_{\ga}$ consisting of commuting roots
such that $N$ is $\Gamma_{\ga}'$-injective. Since $S$ is
$\ga$-bounded, by the Lie algebra version of the statement
(Proposition 4.8 in \cite{M}), $R:=\Psi_{\Gamma_{\ga}}^\nu (S)$ is
a simple $\ga$-module for some $\nu \in \gh^*$. Then by Lemma
\ref{commut}, $L_{\gp}(R) \simeq \Psi_{\Gamma_{\ga}}^\nu (M)$.  On
the other hand ${\cal E}_{\ga}(S)[t_0]$ contains a
$\Gamma_{\ga}'$-injective simple submodule $N_0'$, and in
particular $\Psi_{\Gamma_{\ga}'}^\nu (N_0') \simeq R$. Therefore,
if $N':= L_{\gp}(N_0')$, part (i) of the theorem implies ${\cal
E}_{\ga}(N') \simeq {\cal E}_{\ga}(M) $. If $r_0 = t_0$ we  may
assume that $N' = N$ which implies the statement. Suppose now $r_0
\neq t_0$.  Then $N$ and $N'$ are nonisomorphic $\gg$-submodules
of ${\cal E}_{\ga}(M)$, and thus $ \deg_{\ga} {\cal
E}_{\ga}(M)[r_0] \geq \deg_{\ga} N[r_0] +  \deg_{\ga} N'[r_0] >
\deg_{\ga} {\cal E}_{\ga}(N')[r_0] = \deg_{\ga} {\cal
E}_{\ga}(M)[r_0]$, which is impossible.

(iii) Theorem 13.3 in \cite{M} implies that there exists $\lambda \in \gh^*$ such that  $ \Psi_{\Gamma_{\ga}}^\nu (L_{D_{\ga}}(\lambda)) \simeq S$. Now using Lemma \ref{commut} we complete the proof.
\hfill $\square$

\bigskip

\noindent {\bf Note.} The above theorem has been proved in
\cite{M} for the case when $\gg = \ga$. The case when $\gg$ is a
Lie superalgebra of type I and $\ga = \gg_0$ is considered in
\cite{G}. One should note that  the statement of Theorem 8 in
\cite{G} requires a slight modification, namely that the module
$M$  should be assumed to be strictly bounded instead of
infinite dimensional bounded (this modification is necessary
for the cases $\gg = \gg \gl (m|n)$ and $\gg = \gs \gl (m|n)$ and
$m,n \geq 2$).

\bigskip
Recall that the Harish-Chandra homomorphism embeds the center $Z(\gg)$ of $U$
into the $W$-invariant polynomials of $\gh^*$. As usual, this
enables us to assign to each weight $\mu$ a central character
$\theta^{\mu}$, i.e. a homomorphism $\theta^{\mu} : Z(\gg) \to K$.

\smallskip

\noindent {\bf Remark.} Let $\gg$ be one of the following $\gs \gl (m|p)$, $\gp \gs \gl (m|m)$, 
$\go \gs \gp (2|p)$, i.e. $Z(\gg) \neq K$.
One easily checks that the map $\Psi_{\Sigma_{\ga}}^{\nu}$
preserves central characters (in the strong sense). Therefore all
submodules of ${\cal E}_{\ga}(M)$ have the same central character.
As a corollary of Proposition \ref{cohext}, (iii), we see that the
category ${\cal WT}^{\theta^{\lambda}}[t]$ of weight modules  with
fixed central character $\theta^{\lambda}$ whose support lie in a
fixed $Q$-coset $t$ contains finitely many simple objects. This is
true because  the corresponding category ${\cal
O}^{\theta^{\lambda}}[t]$ whose simple objects are highest weight
modules has the same property. In particular, if $B_{\ga}'$ is
fixed basis of $\Delta_{\ga}$ and ${\cal K}^{\lambda}$ is the set
of all couples $(B', \lambda')$, where $B'$ is a basis of $\Delta$
containing $B_{\ga}'$ and $\lambda'$ is a weight for which
$\theta^{\lambda} = \theta^{\lambda'}$, then the ``universal'' semisimple
$(\gg, \ga)$-coherent family ${\cal U}(\theta^{\lambda}, \ga):=\oplus_{(B', \lambda') \in {\cal
K}^{\lambda}} {\cal E}_{\ga}(L_{B'}(\lambda'))$ contains all
simple weight $\gg$-modules with central character
$\theta^{\lambda}$ which are $\ga$-torsion-free. Finally, for $\lambda \in \gh^*$ denote by $SA(\lambda)$ the set of all reductive Lie algebras $\gc$ of $\gg_0$ containing $\gh$ for which $L_{B'}(\lambda)$ is strictly $\gc$-bounded for some basis $B'$ of $\Delta$. Then the module ${\cal U}(\theta^{\lambda}):= \oplus_{\gc \in SA(\lambda)} {\cal U}(\theta^{\lambda}, \gc)$ is a finite direct sum of coherent families and contains all simple weight modules with central character $\theta^{\lambda}$.

\smallskip

\section{On the Composition Series of Weight Modules}

In this section we find a necessary and sufficient condition for
the modules $L_{\gp}(S)$ and $M_{\gp}(S)$ to be  $\gd$-semisimple
(and in particular $\gg_0$-semisimple) provided that the answer
for this  question is known for the 
modules $L_{B}(\lambda)$ and $M_{\gp}(\lambda)$. An interesting
application is obtained for the case $\gg = \go \gs \gp (2|n)$. We
keep the notations and conventions from the previous sections, and
unless otherwise stated we set $D:=B$ in part (iii) of Proposition
\ref{cohext}.

\begin{lemma} \label{irrhw}
Let  $\alpha \in \Delta$ and $\lambda$ and $\lambda + \alpha$ be  $\Gamma_{\ga}$-injective weights. Then for $\nu \in \gh^*$ the following are equivalent:

(i) $\Psi_{\Gamma}^{\nu}(L_{B_{\ga}}(\lambda))$ is  simple;

(ii) $\Psi_{\Gamma}^{\nu}(L_{B_{\ga}}({\lambda}+ \alpha))$ is simple;

(iii)$\Psi_{\Gamma}^{\nu}(L_{B_{\ga}}(\lambda))$ is torsion-free;

(iv) $\Psi_{\Gamma}^{\nu}(L_{B_{\ga}}({\lambda}+ \alpha))$ is torsion-free.
\end{lemma}
\noindent {\bf Proof.} A Theorem of Fernando implies that a simple $\ga$-module $R$ is torsion-free if and only if $\supp(R) = \lambda + Q_{\ga}$ for some $\lambda \in \gh^*$ (see Corollary 1.4 in \cite{M}). Therefore $(i)$ and $(ii)$ are equivalent to $(iii)$ and $(iv)$ respectively. We show now that $(i)$ and $(ii)$ are equivalent. 
Following \cite{M}, for $\mu \in
\gh_{\gs}^*$, denote by $Sing(\theta^{\mu})$ the subset of
$T_{\gs}^*$ consisted of  the cosets $\nu + Q_{\gg}$ for which
$\Psi_{\Gamma}^{\nu}(L_{B_{\ga}}(\mu))$ is not simple. Mathieu
showed that $Sing(\theta^{\mu})$ depends only on the central
character $\theta^{\mu}$ and equals the finite union $\cup_{s \in
W_{\gs}, i \in {\cal F}}(s(\mu) + {\cal H}_i)$ of hyperplanes in
$T_{\gs}^*$. Each hyperplane $s(\mu) + {\cal H}_i$ contains
$s(\mu)$ for some $s \in W_{\gs}$ (if $\mu$ is integral then all
hyperplanes contain $\mu$). The index set ${\cal F}$ is finite and
for its precise definition we refer the reader to  Section 10 in
\cite{M}. We must show that $\lambda^{\gs} + \alpha^{\gs} +
\nu^{\gs} \notin Sing(\theta^{\lambda^{\gs} + \alpha^{\gs} })$
whenever  $\lambda^{\gs} + \nu^{\gs} \notin
Sing(\theta^{\lambda^{\gs} })$. But for $s \in W_{\gs}$ we have
 $s(\alpha^{\gs}) - \alpha^{\gs} = s(\alpha) - \alpha \in
Q_{\gs}$ as $s(\alpha^{\gz}) = \alpha^{\gz}$. Therefore
$\lambda^{\gs} + \alpha^{\gs} + \nu^{\gs} + Q_{\gs} \in
Sing(\theta^{\lambda^{\gs} + \alpha^{\gs} })$ if and only if
$\lambda^{\gs} + \alpha^{\gs} + \nu^{\gs} - s(\lambda^{\gs} +
\alpha^{\gs}) + Q_{\gs} \in \cup_{i \in {\cal F}}{\cal H}_i$ for
some $s \in W_{\gs}$. The latter is equivalent to $\lambda^{\gs} +
\alpha^{\gs} \in  Sing(\theta^{\lambda^{\gs}})$ which completes
the proof. \hfill $\square$

\begin{proposition} \label{JH}
Let $\lambda$ be a  $\Gamma_{\ga}$-injective weight for which
$\Psi_{\Gamma}^{\nu}(L_{B_{\ga}}({\lambda}))  \simeq S$ is a
torsion-free $\ga$-module for some $\nu \in \gh^*$. Let $M$
be one of the $\gd$-modules $L_{B}(\lambda)[t]$ or
$M_{\gp}(\lambda)[t]$, where $t \in T_{\gd}^*$. Then if $0 = M_1^1
\subset ... \subset M_{r_1}^1 \subset M_1^2 \subset ... \subset
M_{r_2}^2 \subset ... \subset M_1^t \subset ... \subset M_{r_t}^t
= M$ is a $\Gamma_{\ga}$-injective composition series of $M$,
$\Psi_{\Sigma}^\nu (M_1^1) \subset \Psi_{\Sigma}^\nu(M_1^2)
\subset ... \subset \Psi_{\Sigma}^\nu(M_1^t) =
\Psi_{\Gamma}^{\nu}(M) $ is a  composition series of
$\Psi_{\Gamma}^{\nu}(M)$ (i.e. of $L_{\gp}(S)[t]$  or
$M_{\gp}(S)[t]$).

\end{proposition}
\noindent
{\bf Proof.} The proposition follows directly from Lemmas \ref{commut} and \ref{irrhw}, and Theorem \ref{Psi}.
\hfill $\square$

If $K$ and $L$ are weight $\gd$-modules, then $K$ is a $\ga${\it
-partially finite enlargement} of $L$ if there exists an
exact sequence $0 \to L \to K \to L' \to 0$ for some
$\ga$-partially finite $\gd$-module $L'$.

\begin{corollary} \label{s-s}
The weight $\gd$-module $L_{\gp}(S)$ (respectively, $M_{\gp}(S)$)
is semisimple if and only if $L_B(\lambda)$ (resp.,
$M_{\gp}(\lambda)$) is an $\ga$-partially finite enlargement of a
semisimple $\gd$-module.
\end{corollary}
\noindent {\bf Proof.} In order to prove the statement for
$L_{\gp}(S)$, it is enough to consider the module $L_{\gp}(S)[t]$,
where $t \in T_{\gd}^*$. If $L_B(\lambda)[t]/ (\oplus_{i=1}^r
L_i)$ is an $\ga$-partially finite $\gd$-module, then
$L_B(\lambda)[t]$ has a $\Sigma_{\ga}$-injective composition
series (\ref{JH2}) for which $t=r$, $r_1 = ...= r_{k-1} = 1$, and
$M_1^l \simeq \oplus_{i=1}^l L_i$ for $l =1,...,r$. Then combining
Theorem \ref{Psi} and  Lemma \ref{irrhw} we see that the
$\gd$-modules $\Psi_{\Sigma_{\ga}}^\nu(L_i)$ are simple and
$L_{\gp}(S)[t] \simeq \Psi_{\Sigma_{\ga}}^\nu(L(\lambda)[t])
\simeq \oplus_{i=1}^r \Psi_{\Sigma_{\ga}}^\nu(L_i)$. Suppose now
that $L_{\gp}(S)[t] \simeq \oplus_{i=1}^r M_i$ is a semisimple
$\gd$-module, and set $N_i:=\Phi_{\Sigma_{\ga}}^{-\nu} (M_i)$. If
$X_i$ is a simple submodule of $N_i$ then the simplicity of
$M_i$ implies $\Psi_{\Sigma_{\ga}}^\nu (X_i) \simeq M_i \simeq
\Psi_{\Sigma_{\ga}}^\nu (N_i)$. Using  Lemma \ref{inj-fin}
we show that $N_i$ is an $\ga$-partially finite enlargement of
$X_i$. Again by  Lemma \ref{inj-fin}, $L_B(\lambda)[t]$ is an
$\ga$-partially finite enlargement of $\oplus_{i=1}^r (N_i)$ which
completes the proof for $L_{\gp}(S)$. The proof for $M_{\gp}(S)$
is analogous. \hfill $\square$

For a simple $\gg_0'$-module $R$ we set $K(R):=M_{\gg^0 \oplus \gg^1}(R)$ for $\gg \neq ${\bf W}$(n)$, and  $K(R):=M_{W_{\geq 0}}(R)$  (provided that $W_{\geq 1}$ acts trivially on $R$) for  $\gg = ${\bf W}$(n)$. We call $K(R)$ a {\it generalized Kac module}. In the particular case of $S = L_{B_{{\gg}_0'}}(\lambda)$ we obtain the (classical)  {\it Kac module} $K_B(\lambda)$.

\begin{corollary} Let $\gg = \go \gs \gp (2|n)$. Then if $S$ is an $\ga$-strictly bounded simple $\ga$-module,  the $\ga$-module $M_{\gp}(S)$ (and therefore $L_{\gp}(S)$ too) is an $\ga$-partially finite enlargement of a semisimple $\ga$-module. In particular, if $\ga = \gg_0$, every bounded generalized Kac module $K (S)$ (and therefore every bounded simple $\gg$-module) is a partially finite enlargement of a semisimple $\gg_0$-module.

\end{corollary}
\noindent
{\bf Proof.} The torsion-free $\ga$-modules of finite length are semisimple (see Theorem 1 in \cite{BKLM}) and therefore the statement is true for a torsion-free $\ga$-module $S$. By applying twice Corollary \ref{s-s} for $\gd = \ga$, we prove the statement for any $\ga$-strictly bounded module $S$.
\hfill $\square$

\section{Explicit Character Formula}
In this section we apply  Theorem \ref{Psi} to the case when $\gd
= \gg$ and $M = M_{\gp}(\lambda)$, and write down an explicit
character formula for any simple weight $\gg$-module $L_{\gp}(S)$
in terms of the multiplicities $[M_B(\nu):L_B(\mu)]$.

\begin{lemma} \label{mblb}
 If $\mu \in \gh^*$ then  $$ [M_B(\lambda): L_B(\mu)] = \sum_{\nu} [M_{B_{\ga}}(\lambda): L_{B_{\ga}}(\nu)] [M_{\gp}(\nu): L_B(\mu)].$$
\end{lemma}
\noindent
{\bf Proof.} Let $$0 = M_0 \subset M_1  \subset...\subset M_k \subset M_{k+1} = M_{B_{\ga}}(\lambda)$$ be a composition series of the $\ga$-module $M_{B_{\ga}}(\lambda)$ with successive simple quotients $M_{i}/M_{i-1} = L_{B_{\ga}}(\mu_i)$ for some $\mu_i \in \gh^*$. Since $M_{\gp}$ is an exact functor, the filtration  $$0 = M_{\gp} (M_0) \subset M_{\gp} (M_1)  \subset...\subset M_{\gp} (M_k) \subset M_{\gp}(M_{k+1}) = M_{B}(\lambda)$$ has successive quotients $M_{\gp}(\mu_i)$. Therefore for every $\lambda, \nu \in \gh^*$, $[M_B(\lambda): L_B(\mu)] = \sum_{\nu} (M_B(\lambda): M_{\gp}(\nu))[M_{\gp}(\nu): L_B(\mu)]$, where $(M_B(\lambda): M_{\gp}(\nu))$ stands for the number of occurencies of $M_{\gp}(\nu)$ as a quotient in the above filtration. Using the exactness of $M_{\gp}$ again we see that $ (M_B(\lambda): M_{\gp}(\nu)) = [M_{B_{\ga}}(\lambda): L_{B_{\ga}}(\nu)]$, which completes the proof.
\hfill $\square$

\medskip

We are now ready to write down a character formula for an arbitrary simple $\gg$-module 
$M \simeq L_{\gp}(S)$. Let us first provide a parametric description of $M$. Proposition \ref{cohext} implies that 
$S \simeq {\cal E}_{\ga}(L_{B_{\ga}}(\lambda))[\sigma + Q_{\ga}]$ for some $\sigma \in \gh^*$. Furthermore, $L_{B_{\ga}}(\lambda) \simeq L_{B_{\ga_1}}(\lambda^{\ga_1}) \otimes ...\otimes  L_{B_{\ga_{k}}}(\lambda^{\ga_k}) \otimes K_{\lambda^{\gz}}$, where $K_{\lambda^{\gz}}$ is the one-dimensional $\gz$-module determined by $\lambda^{\gz}$. The algebras $\ga_i$ are of type $A$ or $C$ (see Proposition \ref{dmp}). We  may choose $B$ so that $B_{\ga_i}$ are the standard bases of $\Delta_{\ga_i}$, i.e. $B_{\ga_i}:= \{ \varepsilon_1 - \varepsilon_2,..., \varepsilon_{m_{i-1}} - \varepsilon_{m_i} \}$ if $\ga_i = \gs \gl (m_i) $ and $B_{\ga_i}:= \{ \varepsilon_1 - \varepsilon_2,..., \varepsilon_{m_{i-1}} - \varepsilon_{m_i}, 2 \varepsilon_{m_i} \}$ if $\ga_i = \gs \gp (2m_i)$. Following sections 8 and 9 in \cite{M} we may choose $\lambda^{\ga_i}$ (and therefore also $\lambda$) as follows. If $\gg = \gs \gl (m_i)$ then $\lambda^{\ga_i} = \sum_{j=1}^{m_i} l_j \varepsilon_j$ where $\sum_j l_j = 0$, and  $l_j - l_{j+1} \in \Z_+ \cup {0}$ if and only if $j = 2,3,...,m_i - 1$. If $\gg = \gs \gp (2m_i)$ then $\lambda^{\ga_i} = \sum_{j=1}^{m_i} l_j \varepsilon_j$ where $l_j \in \frac{1}{2} + \Z$ and $l_1 \geq l_2 \geq ... \geq l_{m_i} \geq -\frac{1}{2}$. The module $M$ is uniquely determined by the ordered tuple $(\gp, \lambda^{\ga_1},...,\lambda^{\ga_k},\lambda^{\gz}, \sigma)$. The elements of this tuple will be the main ingredients in the character formula.

 For a Lie superalgebra $\gk$ with reductive even
part $\gk_0$  and a basis $B_{\gk}$ of $\Delta_{\gk}$ set $a_{\nu,
\mu}^{\gk}:=[M_{B_{\gk}}(\nu): L_{B_{\gk}}(\mu)]$ for any $\nu,
\mu \in \gh_{\gk}^*$. Define $b_{\nu, \mu}^{\gk}$ by the equality
$ \ch L_{B_{\gk}}(\nu) = \sum_{\mu}b_{\nu, \mu}^{\gk} \ch
M_{B_{\gk}}(\mu) $, i.e. $[b_{\nu, \mu}^{\gk}]$ is the matrix
inverse to $[a_{\mu, \nu}^{\gk}]$. For simple Lie algebras $\gk$,
$a_{\nu, \mu}^{\gk}$ and $b_{\nu, \mu}^{\gk}$ are computed using
the Kazhdan-Lusztig algorithm. This was originally conjectured in
\cite{KL}, and proved later in \cite{BB} and \cite{BK}. For $\mu
\in \gh^*$ set $S(\mu):=\Psi_{\Sigma}^{\sigma- \lambda} ( L_{B_{\ga}}(\mu) )$, or equivalently (by Proposition \ref{cohext}) $S(\mu) \simeq {\cal E}_{\ga}(L_{B_{\ga}}(\mu))[\sigma + Q_{\ga}]$.
In particular $S(\lambda) = S$.  We will write the character $\ch L_{\gp}(S)$  in terms of the characters $\ch M_{\gp}(S(\mu))$. In order to compute the character of $M_{\gp}(S(\mu))$ we introduce some notations. 

If $\lambda^{\ga_i} = \sum_{j=1}^{m_i} l_j \varepsilon_j$ we set $\tilde{\lambda}^{\ga_i}:= \sum_{j =1}^{m_i-1}l_{j+1}\varepsilon_j \in  \gh _{\gg \gl (m_i-1)}^*$ for $\ga_i = \gs \gl (m_i)$ and   $\tilde{\lambda}^{\ga_i}:= \sum_{j =1}^{n}(l_{j}+1)\varepsilon_j \in \gh _{\gs \go (2m_i)}^*$ for $\ga_i = \gs \gp (2m_i)$. Note that $\tilde{\lambda}^{\ga_i}$ is dominant integral weight with respect to the standard bases $B_{\gg \gl (m_i-1)} := B_{\gs \gl (m_i-1)}$ and  $B_{\gs \go (2m_i)} := \{ \varepsilon_1 - \varepsilon_2,..., \varepsilon_{m_i -1} - \varepsilon_{m_i}, \varepsilon_{m_i -1} + \varepsilon_{m_i} \}$ of  $\Delta_{\gs \go (2m_i)}$. Now let $\ga_i = \gs \gl (m_i)$. We say that $\lambda^{\ga_i} = \sum_{j=1}^{m_i} l_j \varepsilon_j $ is {\it integral} if $l_1 - l_2 \in \Z$ and {\it nonintegral} otherwise. We call $\lambda^{\ga_i}$ {\it singular} if $l_1 +...+l_j + j = 0$ for some $j$, and {\it regular integral} if it is integral and not singular. For a positive integer $l$ we set  $\mu[l]:= t_{l+1} \varepsilon_1 + t_{1} \varepsilon_2 +...+ t_{l} \varepsilon_{l+1} + t_{l+2} \varepsilon_{l+2}+...+  t_{m_i} \varepsilon_{m_i} - \rho_{B_{\gs \gl (m_i)}}$ whenever $\mu = \sum_{j=1}^{m_i} t_j \varepsilon_j - \rho_{B_{\gs \gl (m_i)}}$, where $\rho_{B_{\gs \gl (m_i)}} = \frac{m_i-1}{2}\varepsilon_1 +...+\frac{1- m_i}{2}\varepsilon_{m_i}$ is the half sum of the positive roots of $\Delta_{\gs \gl (m_i)}$ (we put $\mu[l]:=0$ for $l \geq m_i$). Every regular integral weight  $\lambda^{\ga_i}$ is of the form $\mu[l]$ for some dominant integral weight $\mu$ and $1 \leq l \leq m_i - 1$. We put $d(\lambda^{\ga^i}):= \sum_{j=l}^{m_i -1} \dim (-1)^{j-l}(L_{B_{\gg \gl (m_i-1)}}(\widetilde{\mu [j]}))$ if $\lambda^{\ga_i} = \mu [l]$ is regular integral and $d(\lambda^{\ga_i}):=\dim (L_{B_{\gg \gl (m_i-1)}} (\tilde{\lambda}^{\ga_i}))$ for all other $\lambda^{\ga_i}$. Finally, for $\ga_i = \gs \gp (2m_i)$ set $d(\lambda^{\ga_i}):=\frac{1}{2^{m_i-1}}\dim (L_{B_{\gs \go (2m_i)}} (\tilde{\lambda}^{\ga_i}))$.

The following theorem is our main result.

\begin{theorem} \label{charfor}
We have $\ch L_{\gp}(S) = \sum_{\mu}c_{\lambda, \mu}^{\gg} \ch
M_{\gp}(S(\mu))$, where 
\begin{equation} \label{char}
 c_{\lambda, \mu}^{\gg} = \left\{
\begin{array}{ll}
0  & {\text{if } \mu \text{ is } \ga-\text{partially finite}}, \\
\sum_{\nu = \sum \nu_i} b_{\lambda,\nu}^{\gg} a_{\nu^{\ga_1}, \mu^{\ga_1}}^{\ga_1} a_{\nu^{\ga_2}, \mu^{\ga_2}}^{\ga_2}...
a_{\nu^{\ga_k}, \mu^{\ga_k}}^{\ga_k} & {\text{otherwise}},
\end{array}
\right.
\end{equation}
and $\ch  M_{\gp}(S(\mu)) = \Pi_{\alpha \in \gu^-}(1+e^{\alpha}) d(\lambda^{\ga_1})...d(\lambda^{\ga_k}) \sum_{\beta \in \sigma - \lambda + \mu + Q_{\ga}}e^\beta$.

\end{theorem}
\noindent
{\bf Proof.}  Lemmas \ref{inj-fin} and \ref{mblb},
and Proposition \ref{JH} imply that  $[M_{\gp}(S): L_{\gp}(S(\mu))] = 0$ if $\mu$ is
$\ga$-partially finite, and $[M_{\gp}(S): L_{\gp}(S(\mu))] =
\sum_{\nu} b_{\lambda, \nu}^{\gg} a_{\nu,\mu}^{\ga}$ otherwise. Now observing that $ a_{\nu,\mu}^{\ga} = a_{\mu^{\ga_1}, \eta^{\ga_1}}^{\ga_1} ...a_{\mu^{\ga_k}, \eta^{\ga_k}}^{\ga_k} $ we prove identity \ref{char}. Next, by  Lemma \ref{irrhw} we verify that $S(\mu)$ is a torsion-free module  and by Theorems 11.4 and 12.2 in \cite{M} we have $\deg_{\ga}(S(\mu)) = d(\lambda^{\ga_1})...d(\lambda^{\ga_k})$  which completes the proof.   \hfill $\square$

\medskip
\noindent
{\bf Remark.}  The sum in (\ref{char}) is locally finite, i.e. $\dim L_{\gp}(S)^{\eta}$ equals  a sum of finitely many nonzero summands and therefore the formula is well-defined.
\bigskip

In the case when $\gg$ is a Lie algebra, Theorem \ref{charfor} is
a more explicit version of Mathieu's formula (Theorem 13.3 in
\cite{M}) $$\dim L_{\gp}(S)^{\mu} = \Sup_{\zeta \in Q_{\ga}} (\dim
L_{B}(\lambda)^{\sigma - \lambda + \mu + \zeta}),$$ where $\mu \in
\gh^*$.

In general, i.e. for an arbitrary classical Lie superalgebra $\gg$ of
type I or $\gg = ${\bf W}$(n)$, the matrix  $a_{\nu,\mu}^{\gg}$ (and consequently
$b_{\nu,\mu}^{\gg}$) is not known explicitly.  In order to
describe in more details the applications of Theorem \ref{charfor}
for all Lie superalgebras $\gg$ we consider, we need the notions of
typicality and singularity of a weight module, which we recall
briefly below. For more details see \cite{PS1}  and
\cite{PS2}.

A weight $\mu$ is {\it typical} if the Kac module $K_B(\mu)$ is simple. Otherwise it is {\it atypical}. Let $\gb$ be the Borel subsuperalgebra of $\gg$ determined by $B$, and  $\rho_B$ be the half-sum  of the $\gb$-positive even roots and the $\gb$-negative odd roots. If $\mu \in \gh^*$ and $\alpha \in \Delta$, we say that $\mu$ is $(B, \alpha)$-{\it regular} iff $(\mu + \rho_B , \alpha )\neq 0$, where $(\cdot, \cdot )$ is the standard bilinear form on $\gh^*$.  A theorem of Kac implies that for $\gg = \gg \gl (m|n), \gs \gl (m|n), \gp \gs \gl(m|m), \go \gs \gp (2|n)$,  $\mu$ is typical if and only if it is $(B,\alpha)$-regular for all $\alpha \in \Delta_1$ (the original statement of Kac, Proposition 2.9 in \cite{K}, applies for finite dimensional modules, but one extends his proof easily for arbitrary highest weight modules). Therefore $\mu$ is atypical iff it is $(B,\alpha)$-{\it atypical} for some $\alpha \in \Delta_1$, i.e. $(\mu + \rho_B , \alpha ) = 0$. In the cases  $ \gg = \pp(m)$ or  $\gg = \gs \gp(m)$ Corollary 5.8 in \cite{S1} implies that $\mu$ is atypical if and only if $\Pi_{i\neq j}((\mu + \rho_{B_{\gs \gl (m)}}, \varepsilon_i - \varepsilon_j)-1) = 0$. Finally, for $\gg =${\bf W}$(n)$, by a result of Serganova, \cite{S3}, we have that $\mu$ is atypical if and only if $\mu = a\varepsilon_i + \varepsilon_{i+1} +...+ \varepsilon_{n}$ for some $a \in k$. The weight $\mu$ is {\it singular} if $s(\mu + \rho_B) = \mu + \rho_B$ for some $s \in W$. Note that in the case when $Z(\gg) \neq K$, i.e. for all Lie superalgebras except $\gs \pp (m)$, $\pp(m)$ and {\bf W}$(n)$, we define $\theta^{\mu}: Z(\gg) \to K$ to be 
{\it typical central character} (respectively, {\it atypical} or {\it singular central character}) if $\mu$ is typical (resp., atypical or singular) weight. If $Z(\gg) \neq K$, a simple $\gg$-module with central character $\theta^{\mu}$ is  {\it typical module} (respectively, {\it atypical} or {\it singular module}) if $\theta^{\mu}$ is typical (resp., atypical or singular).

It is clear that for a typical weight $\nu$ we have $a_{\nu,\mu}^{\gg} = a_{\nu,\mu}^{\gg_0}$, in which case Theorem \ref{charfor} provides a character formula for $M$. The applications of Theorem \ref{charfor} for the atypical case are listed below.

Let $\gg = \gs \gl (m|n)$ or $\gg = \gg \gl (m|n)$. In this case there is a conjectured formula due to J. Brundan, written in terms of canonical bases of the quantized universal enveloping algebra $U_{q}(\gg \gl (\infty))$, see \cite{B}. Therefore, our result reduces the character problem for $\gs \gl (m|n)$ to the Brundan conjecture. In the case when $n=1$, i.e. for $\gg = \gs \gl (m|1)$ the matrix $a_{\nu,\mu}^{\gg}$ has been found by Serganova (Corollary 9.3 in \cite{S1}) for a set of weights $\nu$ which includes all atypical nonsingular weights.  More explicitly,  if $\nu$ is a nonsingular $\bar{\alpha}$-atypical and $\bar{\alpha} \in B \cap \Delta_1$, then $a_{\nu,\mu}^{\gs \gl (n|1)} = a_{\nu, \mu} ^{\gg \gl (n)} + a_{\nu - \bar{\alpha}, \mu} ^{\gg \gl (n)}$ which provides a character formula for all simple nonsingular atypical $\gs \gl (m|1)$-modules.

Let $\gg = \gp \gs \gl (m|m)$. The question of finding $a_{\nu,\mu}^{\gg}$ for atypical weights $\nu$ is open in general. 

Let $\gg = ${\bf p}$(m)$ or $\gg = \gs${\bf p}$(m)$.   In this case also we don't know much about atypical highest weight modules. We obtain another version of  the character formula  established by Serganova (Corollary 5.8 in \cite{S2}).

Let $\gg = \go \gs \gp (2|n)$, $n=2q$. Serganova's reduction method for atypical representations applies for $\gg$ (Section 4 in \cite {N}). As in the case of $\gg = \gs \gl (m|1)$, this method leads to a character formula for all nonsingular highest weight modules,  but  the formula is not available in the literature yet (\cite {NS}). Then Theorem \ref{charfor} yields a character formula for all nonsingular atypical simple $\go \gs \gp (2|n)$-modules.

Let $\gg =${\bf W}$(n)$. For simplicity we write $W$ instead of {\bf W}$(n)$.  In this case we have a character formula for all simple $W$-modules. A weight $\lambda \in \gh_W^*$ is atypical if and only if  $\lambda = a\varepsilon_i + \varepsilon_{i+1} +...+ \varepsilon_{n}$ for some $a \in K$ and $1 \leq i \leq n$, see \cite{S3}. Fix  $\lambda = a\varepsilon_i + \varepsilon_{i+1} +...+ \varepsilon_{n}$ and set:
$$ s_{\lambda, \mu}^{W} := \left\{
\begin{array}{ll}
\sum_{j\geq 0} (-1)^j \delta_{\lambda - j \varepsilon_i, \mu},  & a \notin \Z_{+}, \\
\sum_{j=0}^{a-1} (-1)^j \delta_{\lambda - j \varepsilon_i, \mu} - \sum_{k\geq 0}(-1)^{a+k} \delta_{- k \varepsilon_n, \mu} + \sum_{l>0} (-1)^{a+l} \delta_{\lambda - l\varepsilon_{i-1} - (a-1)\varepsilon_i,\mu}, & a \in \Z_{+},
\end{array}
\right.$$
\noindent
where $\delta_{\nu,\eta} = 0$ if $\nu \neq \eta$ and $\delta_{\nu,\eta} = 1$ if $\nu = \eta$. A result of Serganova, \cite{S3}, implies that $\ch L_B(\lambda) = \sum_{\mu} s_{\lambda,\mu}^{W} \ch K_B(\mu)$. Now applying Lemma \ref{mblb} for $\gp:= W_{\geq 0}$ we verify that 
$b_{\lambda, \mu}^{W} = \sum_{\nu} s_{\lambda, \nu}^{W} b_{\nu, \mu}^{\gg \gl (n)}$ and finally obtain the character formula for $M$ by Theorem \ref{charfor}.

\noindent
Author's address:\\
Department of Mathematical and Statistical Sciences\\
University of Alberta \\
Edmonton, Alberta T6G 2G1 \\
CANADA\\
\medskip
e-mail address: grant@math.ualberta.ca

\end{document}